# A non-iterative method for optimizing linear functions in convex linearly constrained spaces

**GERARDO FEBRES**[*,**]

[*] *Departamento de Procesos y Sistemas, Universidad Simón Bolívar, Venezuela.*
[**] *Laboratorio de Evolución, Universidad Simón Bolívar, Venezuela.*

*This document introduces a new strategy to solve linear optimization problems. The strategy is based on the bounding condition each constraint produces on each one of the problem's dimension. The solution of a linear optimization problem is located at the intersection of the constraints defining the extreme vertex. By identifying the constraints that limit the growth of the objective function value, we formulate linear equations system leading to the optimization problem's solution. The most complex operation of the algorithm is the inversion of a matrix sized by the number of dimensions of the problem. Therefore, the algorithm's complexity is comparable to the corresponding to the classical Simplex method and the more recently developed Linear Programming algorithms. However, the algorithm offers the advantage of being non-iterative.*

## 1. INTRODUCTION

Since the apparition of Dantzig's Simplex Algorithm in 1947, linear optimization has become a widely used method to model multidimensional decision problems. Appearing even before the establishment of digital computers, the first version of the Simplex algorithm has remained for long time as the most effective way to solve large scale linear problems; an interesting fact perhaps explained by the Dantzig [1] in a report written for Stanford University in 1987, where he indicated that even though most large scale problems are formulated with sparse matrices, the inverse of these matrixes are generally dense [1], thus producing an important burden to record and track these matrices as the algorithm advances toward the problem solution.

Other algorithms have been implemented radically different strategies from the simplex method. Some of them deserve to be mentioned. In 1965 J. Nelder and R. Mead [2] introduced an algorithm that evaluates a growing simplex's vertexes to find the optimal value of minimization problems. This algorithm disregards derivatives of the objective function and thus is affected when the problem's dimension grows. During the 70's Shamos and Hoey [3,4] and Shamos [5] explored options for building efficient linear optimization algorithms based on the geometrical properties of convex spaces. Their proposal mostly included problems of two and three dimensions.

The first polynomial time algorithm for solving an LP was announced by Khanchian [6] in 1979. The algorithm is based on building a hyper-ellipsis, which volume and shape evolve in order to contain an extreme feasible point. Under certain configurations of the constraints, the ellipsoidal methods' convergence process may become exponential. This probably is one of the few weaknesses of these methods.

The next class of algorithms for solving LPs are called 'Interior Point' methods. As the name suggests, these algorithms start by finding an interior point of the constraint-polytope and then proceeds to the optimal solution by moving inside the polytope. The first interior point method was presented by Karmarkar [7] in 1984. This work unchained an important debate about the patentability of math formulas and software. Karmarkar was granted a patent for his algorithm, but the patenting of software and ideas is not allowed any more in the United States. Karmarkar's Patent expired in 2006 and now his algorithm is in the public domain. More recently a method presented by Kelner and Spielman in 2006 [8] runs in polynomial time for all inputs.

Another branch of development started with the named Las Vegas Algorithm introduced by K. L. Clarkson in 1995 [9]. Clarkson's algorithm recursively identifies redundant constraints and afterwards determines the optimal point by intersecting the remaining non-redundant constraints. Due to the recursive process, the algorithm is effective for problems with a small number of dimensions. Later modifications improve the algorithm's performance as shown in [10]. Besides the methods for solving Linear Optimization Problems, some techniques have been developed to improve their performance. We first mention is the Scaling of Linear Problems is used prior to the optimization process itself. This technique deals with the constraint coefficients in the base matrix to reduce the number of iterations required to solve linear optimization problems. Even though in general, these methods actually reduce the number of iterations, their effectiveness is a matter of discussion because often they add more computational costs than the complexity reduction they offer. See Elble J. M, and Sahinidis N. V. [11] for a discussion on this point. Improving the phases of the Simplex algorithm by limiting the complexity of dealing with large number of constraints continues to be a subject of study. T. D. Hansen and U. Zwick [12] offer a method to improve the pivoting process applied in the Simplex algorithm.

Despite the advantages some of these improvements and modifications may offer, the Simplex method remains as the most widely used linear programming algorithm, and the disadvantage of being an iterative algorithm has not been a reason to substitute it. Yet, we believe a non-iterative algorithm would be convenient. This paper presents a different method for solving linear optimization problems. The strategy is based on some characteristics of the boundary on any convex multidimensional polytope. Our major objective is to develop a non-iterative method. An algorithm is devised and presented.

## 2. METHODS

The method's strategy is to deal with the complexity of an n-dimensional optimization problem by inspecting the objective function's bounding conditions in the sense of every single orthogonal dimension. Two different classes of constrains are defined. All constraints are classified to ease implementing the method. After performing geometrical computations we are able to recognize the active constraints. The extreme-vertex coordinates are determined by solving a system of linear equations formulated with the recognized the active constraints.

### 2.1 The general approach and definitions

Consider the classical linear optimization problem:

$$\begin{aligned}&\max\quad \mathbf{c}^T \cdot \mathbf{x} &(1)\\ &\text{subject to:}\quad \mathbf{A}\,\mathbf{x} \leq \mathbf{b}\end{aligned}$$

where the vector $\mathbf{c} \in \mathbb{R}^n$ describes the objective function's coefficients, and $\mathbf{b} \in \mathbb{R}^m$ and $\mathbf{A} \in \mathbb{R}^m \times \mathbb{R}^n$ are the resource vector and the constraint coefficients which describe the feasible region. Vector $\mathbf{x} \in \mathbb{R}^n$ is the variable vector where the optimization problem is defined.

A way to identify the extreme vertex is to inspect the gradients of the objective function and the constraints. We cannot compute the gradient at the vertex formed by the intersection of several surfaces. However, we do know that when taking director vector from one of the hyperplanes surrounding the vertex to another, the vector's direction suddenly 'jumps' from one to another as it passes over the vertex where the hyperplanes intersect. Thus we can imagine an *n*-dimensional solid-angle forming with the director-vectors associated to the *n* hyperplanes intersecting at the vertex. We refer to these angles as the 'transition-angles'. Figure 1 illustrates vertexes in spaces of two and three dimensions and the corresponding director-vectors-formed two and three dimensional angles.

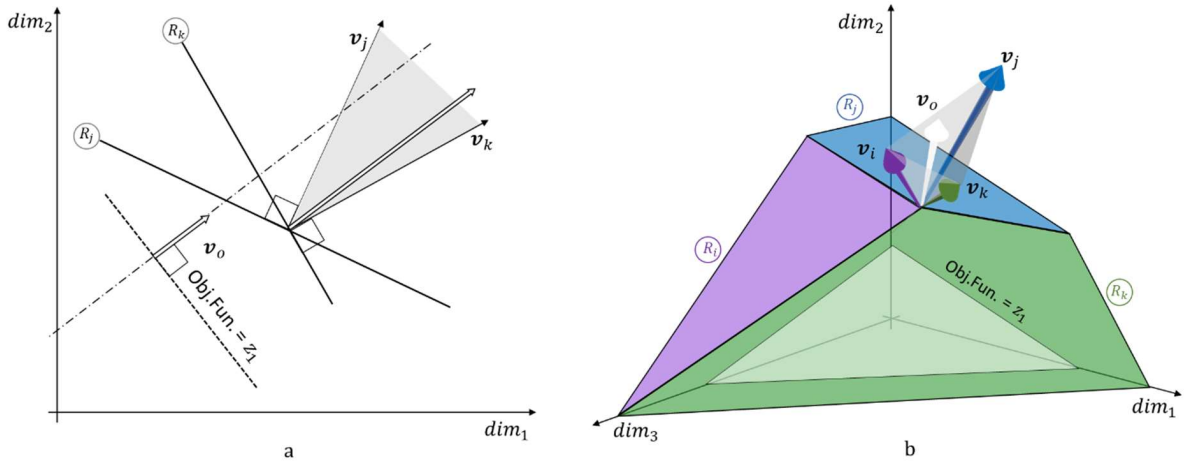

**Figure 1.** Formation of transition–angles associated to the optimal vertex in two and three dimensional spaces. In both cases the vector normal to the objective function is contained by the angle formed within the vectors normal to the hyperplanes intersecting at the vertex.

The direction of the vector normal to the polytope at the exact vertex's location is undetermined. But we know that whatever the direction of that vector is, it must be contained in the core of the transition-angle corresponding to the constraints intersecting at the vertex.

**Theorem 1. Unicity of the constraint-formed hyper-solid-angle containing the objective-function's directive vector.** In a *n*-dimension linear maximization problem with an objective function and a convex closed feasible space, the objective function's normal vector passes through and only passes through (is contained in) the hyper-solid-angle defined by normal vectors to the hyperplanes intersecting at the extreme vertex, upon the objective-function's maximum value.

**Proof:** The Fundamental Theorem of Linear Programming directly proves the existence of a vertex representing the location of the optimal solution for a linear optimization problem. Since this vertex is located where a linear objective function is maximum for polytope that is convex, it must not be more than one vertex. Being this an extreme point (from where the objective-function's value cannot further grow), the objective-function's normal vector must be oriented within the hyper-solid angle formed by the normal vector of the constraints intersecting the extreme vertex. □

The problem solution strategy is then to identify the constraints intersecting at the vertex where the maximum objective function value is. A possible way to identify these constraints is by considering the distances from the constraints to the origin; the farthest would be candidates to participate in the definition of the extreme vertex. However, this may not be a direct identification because some constraints, especially those being far from the origin, tend to be redundant constraints that do not define in any way the feasible region. To handle this situation we proceed to organize the problem formulation and setting some definitions regarding a convenient way of formulating the problem.

**Definition: Inward constraint.** In an $n$-dimensional space, a constraint represented by an inequality involving a $(n-1)$-dimensional surface that splits the $n$-space into two subspaces, is said to be an inward constraint if the vector normal to the constraint's boundary and pointing toward the valid sub-space, forms an angle with the objective-function's normal vector that is greater than $\pi/2$.

**Definition: Outward constraint.** In an $n$-dimensional space, a constraint represented by an inequality involving a $(n-1)$-dimensional surface that splits the $n$-space into two subspaces, is said to be an inward constraint if the vector normal to the constraint's boundary and pointing toward the valid sub-space, forms an angle with the objective-function's normal vector that is smaller than $\pi/2$. Typical non-negativity constraints are outward constraints for maximization problems.

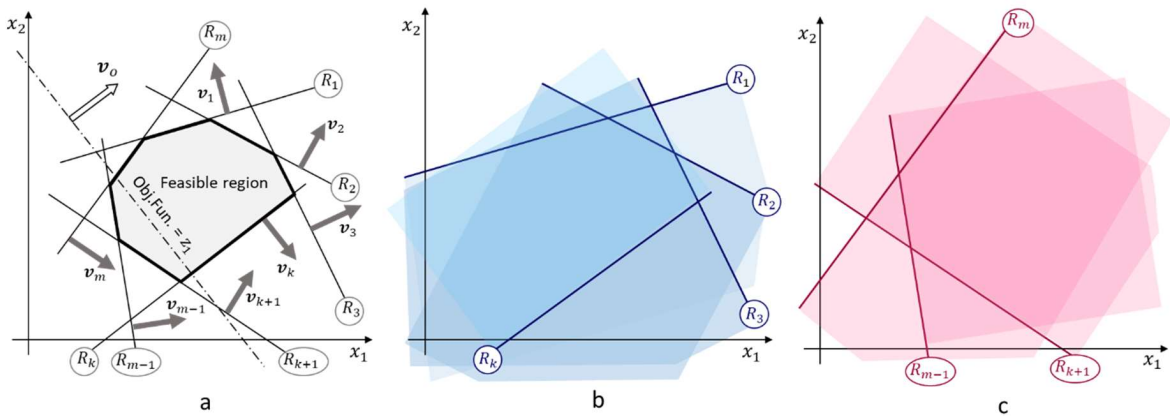

**Figure 2.** Graphic representation of the elements of a generic two-dimensional linear optimization problem. **a**: Constraints and the fastest growing direction of objective function. **b**: Constraints *1* to *k* representing *inward constraints*. **c**: Constraints *k+1* to *m* representing *outward constraints*.

Figure 2 presents a sketch of a two-dimensional optimization problem. Inward and outward constraints are presented in separated graphs. Interpreting Problem (1) geometrically, each constraint represents a hyper-plane, and the set of constraints comprises the envelope of the feasible region forming a polytope in a space of $n$ dimensions. The polytope has vertexes each one located at the intersection of a different group of $n$ constraints. Therefore, there are $\binom{m}{n}$ vertexes in the polytope's surface. This can be a large number of vertexes, but with the exception of some special cases, only one of them represents the conditions where the objective function's value is maximized.

**The inward-outward representation**

After inspecting all constraints, Formulation (1) can be rearranged and written as follows:

$$max \quad z = c_1 x_1 + c_2 x_2 + \ldots + c_n x_n \quad (2)$$

subject to:

$$\sum_{i=1}^{n} a_{ji} x_i \leq b_j \quad (j = 1, \ldots, k) \quad \text{Inward Constraints 1 to k}$$

$$\sum_{i=1}^{n} a_{ji} x_i \geq b_j \quad (j = k+1, \ldots, m) \quad \text{Outward Constraints k+1 to m}$$

Notice this constraint-classification method does not indicate whether or not a constraint is redundant, or on the contrary, makes the problem unfeasible. This classification is useful because typically the optimal point for a maximization problem is defined are inward constraints. Thus, our method primarily searches for inward constraints and then, if necessary, inspects some of the outward constraints for the dimensions not limited by inward constraints.

**Constraints' normal vectors and their angle with the objective function's normal vector**

The unitary normal vector of the objective function is

$$\boldsymbol{v}_o = (v_{o1}, v_{o2}, \ldots, v_{on}) = \frac{1}{\sqrt{\sum_{1}^{n} c_i^2}} (c_1, c_2, \ldots, c_n). \quad (3)$$

We can use the coefficients of constraint $j$ to specify its unitary normal vector. This yields

$$\boldsymbol{v}_j = (v_{j1}, v_{j2}, \ldots, v_{jn}) = \frac{1}{\sqrt{\sum_{1}^{n} a_{ji}^2}} (a_{j1}, a_{j2}, \ldots, a_{jn}). \quad (4)$$

At this point it is worthwhile to emphasize the need for sticking to the format shown in Expression (2), so that the validity the method is ensured.

The cosine of the angle $\alpha_{oj}$ formed by the objective function's director vector $\boldsymbol{v}_o$ and the director vector of constraint $\boldsymbol{v}_j$ can be computed as:

$$\cos \alpha_{oj} = \frac{\boldsymbol{v}_j \cdot \boldsymbol{v}_o}{\|\boldsymbol{v}_j\| \cdot \|\boldsymbol{v}_o\|}. \quad (5)$$

**The By-the-Objective's-Direction Distance (BODD)**

The term by-the-objective's-direction distance (BODD) is used to refer to the distance from the origin to a hyperplane measured over a line that goes along the objective function's normal vector.

Figure 3 illustrates a space where linear constraints named $R_j$ and $R_g$ are represented by solid lines. The dot-dash line shows the maximum growth direction of the objective function passing through the origin. The distances $d_j$ and $d_g$ correspond to the BODDs of constraints $R_j$ and $R_g$ respectively.

## 2.2 The strategy to locate the extreme vertex

Our approach is to solve the linear maximization problem by identifying the constraints defining the extreme vertex. This would be rather easy had it not been for the redundant constraints that may be present in the problem formulation. Detecting redundant constraints is a difficult task. Intuitively results easy to understand that redundant constraints cause great confusion in most algorithms because they intersect other constraints, forming unfeasible vertexes which not being properly dealt, bring the extreme vertex search to a long and fruitless process. In fact, proving or disproving constraint redundancy is computationally equivalent to a linear program [13]. Previous works that may be related to this issues were presented by M. I. Shamos and D. Hoey [3][4]. They classified some optimization problems and created categories they named Closest Point Problems and Geometric Intersection Problems, but in these studies, they do not present problems of more than two dimensions (in the geometrical sense). Later, in his PhD Thesis, Shamos [5] deals with multidimensional spaces, though his treatment focuses on geometric problems and does not directly connect to optimization problems as he had done previously with two dimensions.

The method relies on geometrical properties of convex polytopes and applies them to multidimensional spaces. We begin noting that for a linear maximization problem, the closest constraint with the origin always defines at least a fraction of the boundary of the feasible region. Since the objective function's maximum growth direction is the relevant direction, we select measuring the distance between any constraint and the origin by the direction of maximum growth of the objective function that is the BODD. Therefore, we have a way to identify a constraint, in the objective function's growing direction, that must not be redundant.

**Theorem 2: Non-redundancy of the By-the-objective-direction closest-inward-constraint.** In a linear maximization problem with an extreme vertex exclusively limited by inward constraints, the inward-constraint with the closest by-objective-direction-distance (BODD) is non-redundant.

**Proof:** After adding a constraint that is closer to the origin by the objective function normal direction, the previous BODD closest-constraint may become redundant, but then it will not be the BODD constraint any more. □

**Computing the BODD**

Following the objective-function's growing direction, the coordinates of the points over the objective-function's normal direction are proportional to those described in Equation (3). Therefore, any connection between the origin and a point in the objective-function's normal direction must have the following relationship between its coordinates and the value of coordinate $x_1$:

$$x_2 = \frac{c_2}{c_1} x_1 \ , \ x_3 = \frac{c_3}{c_1} x_1 \ , \ldots, \ x_n = \frac{c_n}{c_1} x_1 \ .$$

The point where objective-function's director vector crosses the constraint $j$ is found joining Expressions (5) and (6). We obtain

$$x_1 = \frac{b_j}{G_j} \ , \quad x_2 = \frac{c_2 \, b_j}{c_1 G_j} \ , \ldots, \ x_n = \frac{c_n \, b_j}{c_1 G_j} \ ,$$

where

$$G_j = a_{j1} + a_{j2}\frac{c_2}{c_1} + a_{j3}\frac{c_3}{c_1} + \cdots + a_{jn}\frac{c_n}{c_1} \ .$$

Finally, the Euclidian distance from the origin and constraint $j$ along the objective-function's director vector is computed as

$$d_j = \frac{b_j}{G_j}\sqrt{1 + \left(\frac{c_2}{c_1}\right)^2 + \left(\frac{c_3}{c_1}\right)^2 + \cdots + \left(\frac{c_n}{c_1}\right)^2} \ . \tag{6}$$

**Computing the coordinates of the by-the-objective-direction-maximum-value-point (BODMP).**

The by-the-objective-function-direction-maximum-value-point (BODMP) is located at the intersection of the objective function's axis and the inward constraint with the minimal BODD. As shown in Figure 3, constraint $R_g$ is the constraint to consider for the BODMP. Scaling vector $\boldsymbol{v}_o$ with Expression (6) for the distance $d_j$ for any constraint $j$ we can determine the coordinates of the BODMP that is over constraint $R_g$. Therefore

$$\boldsymbol{x}_{BODMP} = d_g \cdot (v_{o1}, v_{o2}, \ldots, v_{on}) \ , \tag{7}$$

where we have used the sub-index $g$ to indicate distance $d_g$ refers to the closest BOD distance constraint.

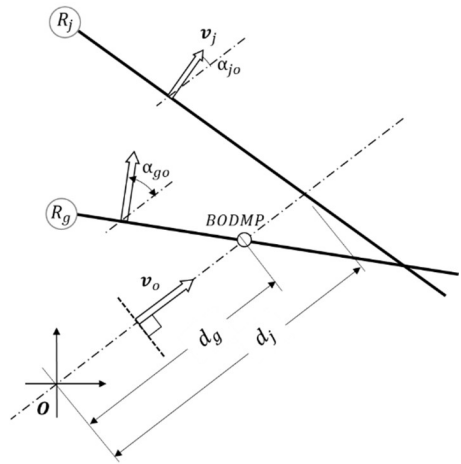

**Figure 3. Identifying the by-the-objective-function-maximal-value-point (BODMP).** Constraints $R_g$ and $R_j$ represented in the context of an $n$-dimensional space with some objective function. The BODMP is located at the intersection of the $\boldsymbol{v}_o$ axis and constraint $R_g$. Constraint $R_g$ defines the BODMP because $d_g < d_j$ among all inward constraints shown.

**Extreme vertex active-constraints determination**

To locate the extreme vertex coordinates, we first look at each one of the problem's dimensions and inspect if the objective function's value can further grow in that dimension or, on the contrary, it is limited by some constraint. In this regard, we state the following theorem:

**Theorem 3: Orthogonal dimension limiting in linear optimization problems.** An $n$-dimensional linear optimization problem is bounded with respect to an objective function if and only if all its orthogonal dimensions are bounded with respect to the same objective function.

**Proof:** If a set of constraints do not bind certain dimension of an $n$-orthogonal-dimensional space, then the indefinite growth of the unbounded dimension can only be stopped by an additional constraint that bounds the unbounded dimension, which now becomes a bounded dimension. Therefore, each one of the $n$ dimensions must be bound in order to have a bounded problem-solution.□

At the BODMP the value associated with the objective function has increased until the closest constraint stops it from growing in the same direction. But then the point evaluated can move in other directions until a constraint limits the growth in the direction of the move. The point evaluated moves while there are still dimensions where additional moves produce an increase in the objective function. This continues until there are no more directions in which the objective function's value can increase. At that point, as Theorem 3 states, there is a limiting constraint indicating the end of the feasible space. This means that determining the set of constraints limiting the problem's dimensions leads us to know the constraints defining the extreme vertex. Once we know the extreme vertex limiting constraints, we can compute the coordinates of the vertex.

An inward-unlimited dimension occurs when the main component of the objective function's vector in that dimension $v_{oi}$, grows faster than the main component $v_{oj}$ of any inward constraint does in the same dimension. Thus, using vector descriptions $\boldsymbol{v}_o$ and $\boldsymbol{v}_j$ in Expressions (3) and (4) we can obtain and compare the vector's main components. Therefore a criterion to test the limiting condition of a dimension is:

$$\text{if for at least one inward constraint } j, v_{oi} \geq v_{ji}, \text{then } dim\ i \text{ is inward} - \text{limited} \tag{8}$$
$$\text{otherwise}, dim\ i \text{ is not inward} - \text{limited}$$

Clearly, when only one inward-constraint limits the feasibility condition in the increasing direction of dimension $i$, this inward limiting constraint is active, and is one of the constraints intersecting at the extreme vertex. When more than one inward-constraint may limit the feasibility condition in the increasing direction of the dimension $i$, we need a criterion to recognize the constraint most likely actually blocks the increase of dimension $i$ values.

Figure 4 shows a plane containing dimension $i$ axis, and the axis containing vector $\boldsymbol{v}_o$ representing the maximum growth direction of the objective function. Inward constraints $R_s$ and $R_l$ are also shown as the inward constraints 'competing' for being the actual limiting element in the dimension $i$ direction. Intuitively, from these constraints, the one closest to the BODMP is the one that most likely actually limits dimension $i$. However, looking for a place with a better objective function's value, the problem's solution intends to move from BODMP to another point by changing the value of coordinate $i$, while values of coordinates in other dimensions that are not specifically shown in the drawing may also change. Thus, other dimensions may be constrained in such a ways that 'force' the change of the prospective point in a specific direction, and encountering an inward constraint that is not the closest to the BODMP, in the Euclidian sense. Nevertheless, this is unlikely to happen, and the loss of precision for this reason should not be major.

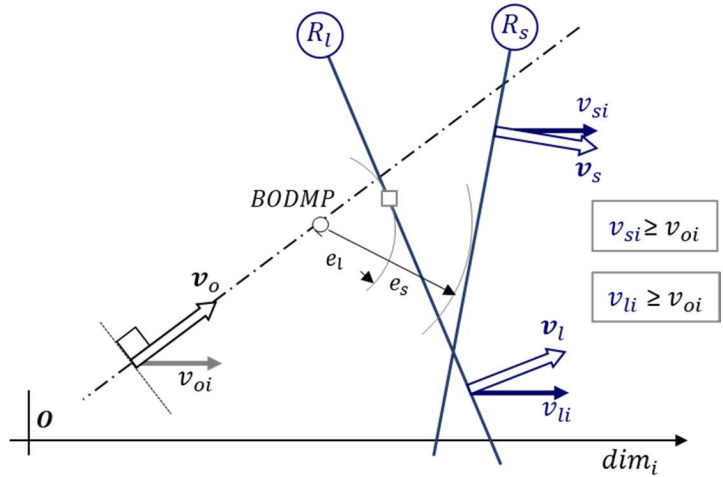

**Figure 4. Inward-constraint limiting condition for dimension $i$.** Constraints $R_l$ and $R_s$ limit the growth of the objective function in dimension $i$. The shorter of the two distances $e_l$ and $e_s$ is the one corresponding to constraint $R_l$, therefore, the active limiting inward-constraint is $R_l$.

Notice the planes containing $R_l$ and $R_s$ are not necessarily perpendicular to the drawing plane. The drawing only shows the lines defined by the intersections of these hyperplanes with the drawing plane. Therefore, we should determine the distances between these hyperplanes and the BODMP, included in Figure 4 as $e_l$ and $e_s$, and shown as if they were contained in the drawing plane; they are not. To determine Euclidian distances from BODMP and constraint $j$ we use

$$e_j = \left| \sum_{i=1}^{n} a_{ji} \cdot (v_{ji} - v_{BODMPi}) \right|. \tag{9}$$

The scalars $v_{BODMPi}$ and $v_{ji}$ correspond to the components of the unit vector pointing to the BODMP and the components of the unit vector normal to constraint $R_j$. The coefficients $a_{ji}$ belong to the description of constraint $R_j$. Finally, when more than one inward-constraints may be the limiting constraint for dimension $i$, we select the one with closest distance $e_j$.

When there is none inward constraint limiting dimension $i$, we proceed to inspect those outward constraints that may limit the growth in the direction of dimension $i$, as shown in Figure 5. The criterion for selecting an outward constraint as the limiting element for this dimension is as follows:

$$\text{if for at least one outward} - \text{constraint } j, v_{oi} \leq v_{ji}, \text{then } dim\ i \text{ is outward} - \text{limited}, \quad (10)$$
$$\text{otherwise,}\ dim\ i \text{ is not outward} - \text{limited}.$$

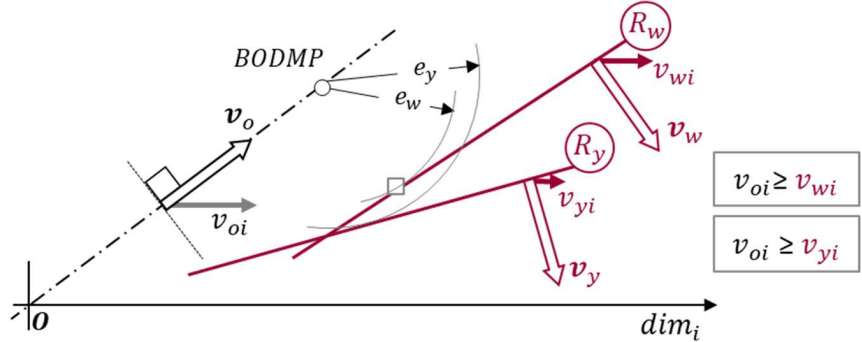

**Figure 5. Outward-constraint limiting condition for dimension $i$.** Constraints $R_w$ and $R_y$ limit the growth of the objective function in dimension $i$. The shorter of the two distances $e_w$ and $e_y$ is the one corresponding to constraint $R_w$, therefore, the active limiting outward-constraint is $R_w$.

## 3.  A PSEUDO-ALGORITHM

An effective strategy to reach the optimal solution for Problem (1) is presented. It is conducted by the following steps:

**Step 0**. Organize the formulation following the format of Problem (2). Classify constraints as inward and outward constraints. Compute the objective function vector $\boldsymbol{v}_o$ unitary component for each dimension. Compute the angles $\alpha_{jo}$ formed by the vector normal to each constraint and the objective function vector $\boldsymbol{v}_o$.

**Step 1**. Considering Theorem 2 and Equation (7) determine the coordinates of the $BODMP$.

**Step 2**. For each dimension recognize the limiting inward-constraint using Criterion (8). If a dimension is not bounded, then recognize the limiting outward constraint using Criterion (10). If for a dimension there is no limiting constraint, then stop; the problem solution is not bounded.

**Step 3**. Build the set of $n$ dimensions that intersect at the extreme point and form matrix $\boldsymbol{B}$ with the coefficients of the extreme-vertex-defining constraints. Also build the resources vector $\boldsymbol{b}$.

**Step 4**. Solve $x^* = \boldsymbol{B}^{-1} \cdot \boldsymbol{b}$.

## 4.  AN EXAMPLE

The following example shows the steps. Consider the 3-dimension problem with 8 constraints:

**Steps 0. Formulation and constraint classification**

| | max | 0.5 | $x_1$ | + | 1.0 | $x_2$ | + | 2.0 | $x_3$ | | | $d_j$ | $\alpha_{jo}$ | Type |
|---|---|---|---|---|---|---|---|---|---|---|---|---|---|---|
| Subject to: | | | | | | | | | | | | | | |
| $R_1$: | | 2.1 | $x_1$ | + | 3 | $x_2$ | + | 1 | $x_3$ | $\leq$ | 5.2 | 1.894 | 0.802 | Inward |
| $R_2$: | | 1.7 | $x_1$ | + | 2.8 | $x_2$ | + | 2.1 | $x_2$ | $\leq$ | 5 | **1.518** | 0.494 | |
| $R_3$: | | 3 | $x_1$ | + | 1 | $x_2$ | + | 2 | $x_2$ | $\leq$ | 5.5 | 1.939 | 0.710 | |
| $R_4$: | | 1.1 | $x_1$ | + | 2.3 | $x_2$ | + | -1 | $x_3$ | $\leq$ | 5.3 | 14.287 | 1.435 | |
| $R_5$: | | 2.1 | $x_1$ | + | 3 | $x_2$ | + | 1.1 | $x_3$ | $\leq$ | 5.8 | 2.126 | 0.776 | |
| $R_6$: | | 1 | $x_1$ | + | 0 | $x_2$ | + | 0 | $x_3$ | $\geq$ | 0 | 0.000 | 1.351 | Outward |
| $R_7$: | | 0 | $x_1$ | + | 1 | $x_2$ | + | 0 | $x_3$ | $\geq$ | 0 | 0.000 | 1.119 | |
| $R_8$: | | 0 | $x_1$ | + | 0.2 | $x_2$ | + | 1 | $x_3$ | $\geq$ | -1 | 4.583 | 0.344 | |

**Step 1. Coordinates of the $BODMP$**

Inward constraint $R_2$ shows the closest distance to the origin (1.518) measured by the objective direction (BOD). The $\boldsymbol{x}_{BODMP}$ coordinates of are computed as follows:

$$\boldsymbol{v}_o = (v_{o1}, v_{o2}, \ldots, v_{on}) = \frac{1}{\sqrt{\sum_1^n c_i^2}} (c_1, c_2, \ldots, c_n) = (0.218, 0.436, \ldots, 0.873).$$

$$\boldsymbol{x}_{BODMP} = d_{BODMP} \cdot (v_{o1}, v_{o2}, \ldots, v_{on}) = 1.518 \cdot (0.218, 0.436, \ldots, 0.873) = (0.331, 0.662, 1.325).$$

**Step 2. By dimension constraint-distance to the BODMP**

| Constraint | Unitary components | | | Distance e | Limiting Constraint | | |
|---|---|---|---|---|---|---|---|
| | $v_{R\,dim1}$ | $v_{R\,dim2}$ | $v_{R\,dim3}$ | To $\boldsymbol{x}_{BODMP}$ | $dim_1$ | $dim_2$ | $dim_3$ |
| Ref. $\boldsymbol{x}_{BODMP}$ | 0.331 | 0.662 | 1.325 | | | | |
| Inward: $R_1$ | 0.553 | 0.790 | 0.263 | 1.781 | Could | Could | Not |
| $R_2$ | 0.437 | 0.720 | 0.540 | 2.576 | Could | Could | Not |
| $R_3$ | 0.802 | 0.267 | 0.535 | **0.564** | **Select** | Not | Not |
| $R_4$ | 0.402 | 0.840 | -0.365 | 0.859 | Could | **Select** | Not |
| $R_5$ | 0.549 | 0.785 | 0.288 | 1.869 | Could | Could | Not |
| Outward. $R_6$ | 1.000 | 0.000 | 0.000 | 0.669 | Not | Redundant | Could |
| $R_7$ | 0.000 | 1.000 | 0.000 | **0.395** | Redundant | Not | **Select** |
| $R_8$ | 0.000 | 0.196 | 0.981 | 0.423 | Redundant | Redundant | Could |

**Step 3. Assemble base matrix $B$**

$$\boldsymbol{B} = \begin{bmatrix} 3 & 1 & 2 \\ 1.1 & 2.3 & -1 \\ 0 & 1 & 0 \end{bmatrix}, \quad \boldsymbol{B}^{-1} = \begin{bmatrix} 0.3743 & -0.5348 & 0.3743 \\ 0.3102 & 0.1283 & -0.0298 \\ -0.7166 & 07380 & 0.0834 \end{bmatrix}, \quad b = \begin{bmatrix} 5 \\ 5.2 \\ 5.5 \end{bmatrix}$$

**Step 4. Solve:**

$$\boldsymbol{x}^* = (x_1^*, x_2^*, x_3^*), \quad \boldsymbol{x}^* = \boldsymbol{B}^{-1} \boldsymbol{b}, \quad \boldsymbol{x}^* = (1.1497, 0.6241, 0.7134),$$
$$z^* = 2.6257.$$

## 5. DISCUSSION

The classical Simplex method strategy is to traverse from a first feasible vertex towards the optimal vertex. The sequence of vertexes visited consists of adjacent vertexes that are conveniently chosen each time a new vertex is selected for evaluation. The number of vertexes may rapidly increase with the number of dimensions and constraints, situating the worst cases of the Simplex algorithm within the category of computational exponential problems. The simplex algorithm indeed visits all $2^n$ vertices in the worst case and this turns out to be true for any deterministic pivot rule [14]. Nevertheless, despite the number of vertexes visited, the relative simple calculations required to update data of neighbor-vertexes at each step, keeps the Simplex method computationally effective for most practical situations, and still, the most popular linear optimization algorithm.

Ellipsoidal methods base their strategy on the behavior of convex polytopes. These methods substitute the traverse of vertexes, with a convergence process that 'directly' builds a hyper-ellipse containing the extreme vertex. This convergence process is polynomial, but may present the problems of stability typical of convergence processes.

The method here proposed is based on characteristics of convex bodies and general geometric properties. Basically, we present a way to find the actual bounding constraints by inspecting a reduced number of them and favorably impacting the resulting complexity of the algorithm. The idea of reducing the original problem's complexity by detecting those constraints not intersecting at the extreme vertex, has been formerly used. In 1995 Clarkson presented the Las Vegas algorithm, based on a random process to reject redundant or non-relevant constraints. There are improvements to this method. One of them can be found in a paper by Brise and Gartner [10]. Yet, these approaches are based on iterative processes. In contrast, we present a non-iterative method to relief the computation by discarding many non-relevant constraints. We believe this is possible thanks to the initial classification of constraints as 'inward' and 'outward', and Theorem 2, that allows the algorithm for a rapid positioning of the pivoting point near the final problem's solution. By decoupling the problem's bounding condition and dealing with simpler bounding problems, one for each dimension, we are able to build a set of active bounding constraints that most likely define the extreme vertex and provides a problem's solution without the need for iterations. Even if some selected constraints are not rigorously the actual limiting constraints, the coordinates of the algorithmic solution will always be close to the coordinates of the real optimum place, as well as the resulting objective function's value, will represent a useful practical solution.

## 6. CONCLUSION

An algorithm to solve optimization linear problems has been presented. The representation of any linear problem as a convex object in an n-dimensions space provides the rational basis for this algorithm. The algorithm's complexity is the same as the complexity of inverting a matrix; $O(n^3)$ or $O(n^{\log_2 7})$ depending on the matrix-inversion algorithm used [15]. In this regard the algorithm does not offer dramatic advantages if compared to the traditional simplex algorithm. However, it does offer the advantage of being a non-iterative method toward the optimal result, which is, therefore, easier to follow and implemented in a wide range of computational environments.

A new strategy for dealing with linear spaces is presented. With this strategy a method for solving linear optimization problems is offered. The method is especially suited for linear optimization problems, but modifications to treat nonlinear objective functions and restrictions are foreseeable. The method relies on rather simple geometrical properties and the intrinsic behavior of convex spaces, thus it may be the basis in the search for strategies to solve optimal linear probabilistic models.